\documentclass[10pt,a4paper]{article}

\usepackage[utf8]{inputenc}
\usepackage[T1]{fontenc}
\usepackage{lmodern}
\usepackage[a4paper,margin=25mm]{geometry}
\usepackage{amsmath,amssymb,amsfonts,mathtools,amsthm}
\usepackage{graphicx}
\usepackage{xcolor}
\usepackage{placeins}
\usepackage{enumitem}
\setlist[itemize]{noitemsep}
\usepackage[numbers,sort&compress]{natbib}
\usepackage[colorlinks=true,allcolors=blue!50!black]{hyperref}
\usepackage{tikz}
\usetikzlibrary{patterns}
\usepackage{pgfplots}
\pgfplotsset{compat=newest}

\newcommand{\logLogSlopeTriangle}[5]
{

    \pgfplotsextra
    {
        \pgfkeysgetvalue{/pgfplots/xmin}{\xmin}
        \pgfkeysgetvalue{/pgfplots/xmax}{\xmax}
        \pgfkeysgetvalue{/pgfplots/ymin}{\ymin}
        \pgfkeysgetvalue{/pgfplots/ymax}{\ymax}

        \pgfmathsetmacro{\xArel}{#1}
        \pgfmathsetmacro{\yArel}{#3}
        \pgfmathsetmacro{\xBrel}{#1-#2}
        \pgfmathsetmacro{\yBrel}{\yArel}
        \pgfmathsetmacro{\xCrel}{\xArel}

        \pgfmathsetmacro{\lnxB}{\xmin*(1-(#1-#2))+\xmax*(#1-#2)} 
        \pgfmathsetmacro{\lnxA}{\xmin*(1-#1)+\xmax*#1} 
        \pgfmathsetmacro{\lnyA}{\ymin*(1-#3)+\ymax*#3} 
        \pgfmathsetmacro{\lnyC}{\lnyA+#4*(\lnxA-\lnxB)}
        \pgfmathsetmacro{\yCrel}{\lnyC-\ymin)/(\ymax-\ymin)} 

        \coordinate (A) at (rel axis cs:\xArel,\yArel);
        \coordinate (B) at (rel axis cs:\xBrel,\yBrel);
        \coordinate (C) at (rel axis cs:\xCrel,\yCrel);

        \draw[#5]   (A)-- node[pos=0.5,anchor=north] {1}
                    (B)--
                    (C)-- node[pos=0.5,anchor=west] {#4}
                    cycle;
    }
}

\title{$\varphi$-FD, a second order finite difference scheme for geometries defined by a level-set function: the Neumann case\thanks{This work was supported by the Agence Nationale de la Recherche, Project PhiFEM, under grant ANR-22-CE46-0003-01, and by the CNRS through the MITI interdisciplinary programs. During the preparation of this work the authors used Claude (Anthropic) to improve the English of the text; they reviewed and edited the content as needed and take full responsibility for it.}}

\author{%
Michel Duprez$^{1,*}$, Vanessa Lleras$^{2}$, Alexei Lozinski$^{3}$, Vincent Vigon$^{4}$, Lisl Weynans$^{5}$\\[10pt]
\begin{minipage}{0.92\textwidth}\small\centering
$^{1}$Université de Strasbourg, CNRS, Inria, ICube, F-67000 Strasbourg, France\\
$^{2}$Université de Montpellier, CNRS UMR 5149, IMAG, 34090 Montpellier, France\\
$^{3}$Université Marie et Louis Pasteur, CNRS, LmB (UMR 6623), F-25000 Besançon, France\\
$^{4}$Université de Strasbourg, CNRS, IRMA (UMR 7501), F-67000 Strasbourg, France\\
$^{5}$Université de Bordeaux, CNRS, Inria, 33405 Talence, France\\[6pt]
\texttt{michel.duprez@inria.fr}, \texttt{vanessa.lleras@umontpellier.fr}, \texttt{alexei.lozinski@umlp.fr},\\
\texttt{vincent.vigon@math.unistra.fr}, \texttt{lisl.weynans@math.u-bordeaux.fr}\\[4pt]
$^{*}$Corresponding author
\end{minipage}
}
\date{}

\begin{document}

\maketitle

\noindent\textbf{Highlights}
\begin{itemize}
\item A new second-order immersed finite-difference scheme for Neumann boundary conditions.
\item No ghost points, cut-cell volumes or interface quadrature; no stabilization term.
\item Two simple precautions make the condition number $O(h^{-2})$, independent of the cut.
\end{itemize}

\section{Introduction}

Finite difference methods on Cartesian grids can approximate partial differential equations on complex geometries when boundary conditions can be imposed without boundary-fitted meshes, typically by embedding the domain into a larger grid and enforcing Neumann conditions through modified stencils, ghost values, or volumetric corrections.
 The classical foundation of this approach is the use of ghost points, in which values outside the physical domain are introduced, and the normal derivatives at the boundaries are approximated using interior finite-difference stencils \cite{leveque}.
 Several authors have proposed related Cartesian-grid techniques for irregular and immersed boundaries. Weynans~\cite{lisl2017} augments Cartesian discretizations with interface unknowns enforcing the boundary conditions, thereby recovering second-order accuracy, consistency, and discrete maximum principles.
 An alternative is to enforce boundary conditions through conservative flux balances on cells cut by the boundary \cite{JOHANSEN199860}. Building on this idea, Arias et al.~\cite{arias2018poisson} introduced a finite volume scheme for Robin (and, in the limit, Neumann) conditions, achieving second-order accuracy for both the solution and its gradient, at the cost of computing cut-cell volumes and fluxes.

In this article, we extend the $\varphi$-FD strategy of~\cite{phiFD}, introduced for the Dirichlet case, which uses a level-set function $\varphi$ to describe the domain and is itself inspired by $\varphi$-FEM~\cite{phifem,neumann,stokes}, to Neumann boundary conditions. The Neumann condition $\partial_n u:=\nabla u\cdot n =g$ on $\partial \Omega$, which reads $\nabla u\cdot\nabla\varphi = g\,\|\nabla\varphi\|$ in terms of the level-set since $n=\nabla\varphi/\|\nabla\varphi\|$, is imposed by relaxing it to $\nabla u\cdot\nabla\varphi-g\,\|\nabla\varphi\|=p\varphi$ near $\partial\Omega$ and eliminating the auxiliary field $p$ by considering this equation at a boundary node and its nearest interior node. The scheme uses only the level-set and interior nodes (no ghost values, flux reconstruction, or cut-cell volumes). The relaxation itself is that of the $\varphi$-FEM Neumann scheme~\cite{neumann}; what is new is its finite-difference form, which needs no integration at all, and the conditioning analysis below. We also identify two elementary precautions --- a normalization of the boundary equations and a rule for the choice of the coupled interior node --- under which the condition number is $O(h^{-2})$ and independent of the smallest cut cell, with no stabilization term.

\section{$\varphi$-FD scheme for Neumann boundary conditions}\label{sec:scheme}

\paragraph{A PDE and a scheme:} In the present article, we consider the Poisson problem with a Neumann boundary condition
\begin{equation}\label{eq:Neumann}
- \Delta u + u = f \text{ in } \Omega\quad\quad 
\partial_n u= g         \text{ on } \partial \Omega,
\end{equation}
  where $\Delta$ is the Laplacian operator, $\partial_n$ the derivative along the outward unit normal $n$ to the boundary,
  $f\in L^2(\Omega)$, $g$ is the trace of a function in $H^1(\Omega)$ and the domain $\Omega$ is described by a level-set function, i.e., $\Omega = \{ \varphi < 0 \}$ so the normal is given by $n=\nabla \varphi /\|\nabla\varphi\|$. For the consistency analysis of Section~\ref{sec:consistency} we assume more: $\varphi$ is $C^3$ near $\partial\Omega$ with $\nabla\varphi\neq0$ on $\partial\Omega$, and $f$, $g$ are smooth enough for $u$ to be $C^4$ up to the boundary (e.g. $f\in C^{2,\beta}(\overline\Omega)$, $g\in C^{3,\beta}(\partial\Omega)$, $\partial\Omega$ of class $C^{4,\beta}$).

Let $h>0$. Define the Cartesian grid $\{x_\alpha = \alpha h \text{ for } \alpha \in\mathbb{Z}^2\}$ and the set of unit directions $\underline{\delta} = \{e_0, e_1, -e_0, -e_1\}$ where $\{e_0, e_1\}$ is the canonical basis of $\mathbb{Z}^2$.
We define the following subsets of the grid:
\[
\Omega_h^{\text{in}} = \{ x_\alpha : \varphi (x_\alpha) < 0 \},~\quad
\partial \Omega_h = \{x_\alpha:x_\alpha\notin \Omega_h^{\text{in}} \text{ but } \exists \delta\in\underline{\delta} : x_{\alpha+\delta} \in \Omega_h^{\text{in}} \} \quad \text{ and }\quad
\Omega_h = \Omega_h^{\text{in}}  \cup \partial \Omega_h.
\]
We denote by $N$ the number of nodes in $\Omega_h$. Using the same notation for functions and their grid evaluation, we define the vectors in $\mathbb R^N$: $f=(f_{\alpha})$, $g=(g_{\alpha})$, $\varphi=(\varphi_\alpha)$, and the vector of unknowns $u=(u_\alpha)$, all indexed by $\Omega_h$. The latter, once determined, will give an approximation of the solution to the Poisson problem. To obtain it, we solve the system of $N$ equations:
\begin{equation} \label{eq:sys}
    \begin{cases}
    -\Delta_hu_\alpha +  u_\alpha    =f_{\alpha} \quad &\forall  x_\alpha \in \Omega_h^{\text{in}}, \\
    \varphi_{\alpha} \big(\nabla_h u_{\alpha_0} \cdot \nabla_h \varphi_{\alpha_0}-G_{\alpha_0} \big) - \varphi_{\alpha_0} \big( \nabla_h u_{\alpha} \cdot \nabla_h \varphi_{\alpha} -G_\alpha\big)= 0\  &\forall  x_\alpha \in \partial\Omega_h,
    \end{cases}
\end{equation}
where
 $\Delta_hu_\alpha$ is the standard 5-point discrete Laplacian
$ (\sum_{\delta \in \underline{\delta}} u_{\alpha+\delta} - 4 u_{\alpha})/h^2$;
 $\alpha_0$ is the index of the nearest node of $\Omega_h^{\text{in}}$ to $x_\alpha \in \partial\Omega_h$ (with the exception stated at the end of this section);
  $G_{\alpha}=g(x_\alpha)\,\|\nabla_h\varphi_{\alpha}\|$ ($g$ being extended to a neighborhood of $\partial\Omega$)
and,
for a vector $v=(v_\alpha)\in \mathbb R^N$, the discrete gradient $\nabla_h v_\alpha$ is the vector $(\partial_{e_0}v_\alpha,\partial_{e_1} v_\alpha)$ where for $\delta\in \{e_0,e_1\}$:
\begin{equation}\label{eq:grad}
  \partial_{\delta} v_{\alpha} = \left\{\begin{array}{ll}
     ({v_{\alpha + \delta} - v_{\alpha - \delta}})/2h, &\text{ if } x_{\alpha + \delta} \in\Omega_h, x_{\alpha - \delta} \in \Omega_h,\\
    ( {- 3 v_{\alpha} + 4 v_{\alpha + \delta} - v_{\alpha + 2 \delta}})/2h, &\text{ if } x_{\alpha +\delta} \in \Omega_h, x_{\alpha + 2 \delta} \in \Omega_h, \text{ but } x_{\alpha-\delta} \not\in \Omega_h,\\
   (  {+ 3 v_{\alpha} - 4 v_{\alpha - \delta} + v_{\alpha - 2 \delta}})/2h, &\text{ if } x_{\alpha -\delta} \in \Omega_h, x_{\alpha - 2 \delta} \in \Omega_h, \text{ but } x_{\alpha+\delta} \not\in \Omega_h,\\
    ( {v_{\alpha_0 + \delta} - v_{\alpha_0 - \delta}})/2h, &\text{ else. }
   \end{array}\right.
\end{equation}
In the last case of \eqref{eq:grad}, we call $x_\alpha$ a boundary-isolated point in the direction $\pm\delta$. Since it does not have enough neighbors in this direction, we use the discrete derivative of its associated interior node $x_{\alpha_0}$. 
We remark that the number of boundary-isolated points is uniformly bounded: such a point needs a grid line meeting $\Omega_h$ in fewer than three consecutive nodes, which can only happen within $O(h)$ of the finitely many points where $\partial\Omega$ is tangent to a grid direction, so their number is bounded independently of $h$. On the domain of Section~\ref{sec:num} there are none on any of the meshes used for the computations; Figure~\ref{fig:results}(a) shows the two that appear on a coarser grid.

\paragraph{Quick justification of our scheme:} We define a function $p$ in a neighborhood of $\partial\Omega$ by the formula:
\begin{equation} \label{eq:scalar_field}
 p(x)\,\varphi(x)= (\nabla u \cdot \nabla \varphi)(x) -G(x)
\end{equation}
with $G(x)= g(x) \|\nabla \varphi(x)\|$.  From classical analysis (explained in the consistency part), the Neumann constraint is equivalent to the continuity of  $p$.  So we can approximate both $p(x_\alpha)$ and $p(x_{\alpha_0})$ by the same value $c_\alpha$. We also replace the gradients $\nabla$ by the discrete gradients $\nabla_h$. From the equation above we derive two approximations: $c_\alpha\varphi_{\alpha_0}=\nabla_h u_{\alpha_0} \cdot \nabla_h \varphi_{\alpha_0} - G_{{\alpha_0}}$ and  $c_\alpha\varphi_{\alpha}= \nabla_h u_{\alpha} \cdot \nabla_h \varphi_{\alpha} - G_{{\alpha}}$. Eliminating $c_\alpha$ gives the second equation of \eqref{eq:sys}.

\paragraph{Two details that govern the conditioning:} Both leave the accuracy of \eqref{eq:sys} untouched but decide whether its condition number depends on the position of $\partial\Omega$ relative to the grid. First, the boundary equations must be normalized individually: their natural magnitude is proportional to $\max(|\varphi_\alpha|,|\varphi_{\alpha_0}|)\,\|\nabla\varphi\|$, so scaling them all by the same power of $h$ gives an artificially small weight to those with $|\varphi_\alpha|,|\varphi_{\alpha_0}|\ll h$; we divide each by its own sup-norm and multiply it by $h^{-2}$, the magnitude of the Laplacian rows, a row scaling that leaves the discrete solution unchanged. Second, $x_{\alpha_0}$ must not lie essentially on $\partial\Omega$: if $\varphi_{\alpha_0}\approx0$ the second term of the boundary equation in \eqref{eq:sys} drops out, which is harmless for a single boundary node but fatal when several share that $x_{\alpha_0}$ --- their equations become proportional and the matrix becomes rank deficient. We therefore take for $x_{\alpha_0}$ the nearest interior node \emph{among those with $|\varphi_{\alpha_0}|\geqslant\tau h$}, $\tau=0.1$, falling back to the nearest one if none qualifies; this keeps $\|x_\alpha-x_{\alpha_0}\|=O(h)$, which is all the consistency analysis uses.

\section{Consistency}\label{sec:consistency}

We start with the homogeneous case $g=0$. Let $\tilde u$ be a smooth extension of the solution $u$ to a neighborhood of $\Omega_h$ whose derivatives up to order 4 are bounded.
For any interior node $x_\alpha \in \Omega_h^{\mathrm{in}}$, 
a standard Taylor expansion yields second-order consistency $-\Delta_h \tilde u_\alpha +  \tilde u_\alpha - f_\alpha = O(h^2)$.

Consider a boundary node $x_\alpha \in \partial\Omega_h$. Suppose first that it is not a boundary-isolated point in either direction, so that $\nabla_h \tilde{u}_{\alpha} = \nabla \tilde u(x_\alpha) + O(h^2)$ and $\nabla_h \varphi_{\alpha} = \nabla \varphi(x_\alpha) + O(h^2)$. Since the homogeneous Neumann condition holds, the function $F(x) := \nabla \tilde u(x) \cdot \nabla \varphi(x)$ vanishes on $\partial\Omega$. 
By definition~\eqref{eq:scalar_field}, $p=F/\varphi$; since $F$ vanishes on $\partial\Omega$ at the same rate as $\varphi$, Hardy's inequality ensures that $p$ is smooth in a neighborhood of $\partial\Omega$ (see e.g.~\cite{phifem}). Hence $F=p\varphi$, with $p(x_\alpha)-p(x_{\alpha_0})=O(h)$ since $\|x_\alpha-x_{\alpha_0}\|=O(h)$.
Since $x_\alpha$ and $x_{\alpha_0}$ are both at distance $O(h)$ from the boundary ($\varphi(x_\alpha),\varphi(x_{\alpha_0})=O(h)$), replacing the discrete gradients by the continuous ones (error $O(h^2)$ at non-isolated nodes) the consistency error of the relaxed Neumann condition decomposes as:
\begin{multline*}
  \varphi_\alpha\,\nabla_h \tilde u_{\alpha_0} \cdot \nabla_h \varphi_{\alpha_0} - \varphi_{\alpha_0}\,\nabla_h \tilde u_\alpha \cdot \nabla_h \varphi_\alpha
  = \varphi(x_\alpha) F(x_{\alpha_0}) - \varphi(x_{\alpha_0}) F(x_\alpha) + O(h^3)\\
  = \varphi(x_\alpha) \varphi(x_{\alpha_0}) \left[ p(x_{\alpha_0}) - p(x_\alpha) \right] + O(h^3) = O(h^3),
\end{multline*}
the $O(h^3)$ coming from the $O(h^2)$ gradient error times $\varphi=O(h)$, the continuous part being exactly $O(h^3)$ by the smoothness of $p$.
Second, if $x_\alpha$ is a boundary-isolated point in one of the directions, the gradient approximation error is $O(h)$, which leads to a consistency error of order $O(h^2)$. Such points are few, however, as remarked after~\eqref{eq:grad}.
Finally, for the non-homogeneous case $g \neq 0$, the same consistency proof holds by defining $F(x) := \nabla \tilde u(x) \cdot \nabla \varphi(x) - G(x)$ with $G(x) = \tilde g(x) \|\nabla \varphi(x)\|$, $\tilde g$ being a smooth extension of $g$, since this function also vanishes on the boundary $\partial\Omega$.

\section{Numerical illustration}\label{sec:num}

\begin{figure}[t]
\centering
\begin{minipage}[c]{0.32\linewidth}\centering
\includegraphics[width=0.95\linewidth]{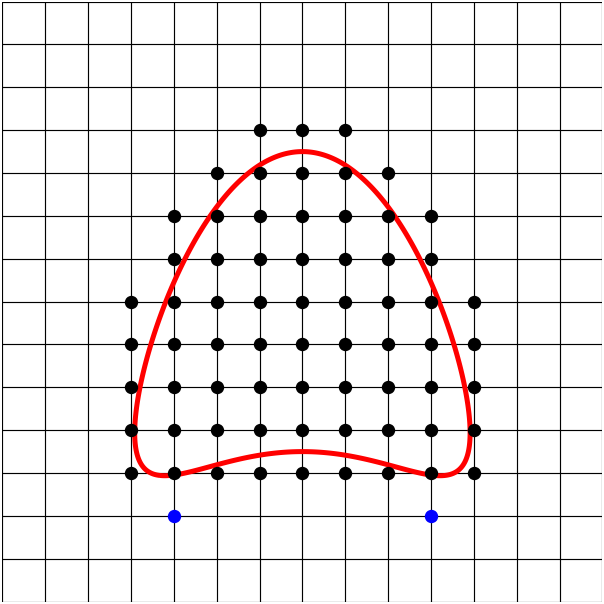}\\[2pt]
{\footnotesize (a) bean-domain $\Omega_h$}
\end{minipage}\hfill
\begin{minipage}[c]{0.62\linewidth}\centering
\begin{tikzpicture}[thick,scale=0.66, every node/.style={scale=1.0}]
\begin{loglogaxis}[xlabel=$h$,xmin=4e-3,xmax=0.15,ymin=4e-5,ymax=0.05,legend style={at={(0.03,0.97)},anchor=north west,font=\scriptsize}]
\addplot[color=blue,mark=triangle*] coordinates {(0.1,3.931633e-02)(0.05,9.680261e-03)(0.025,2.580297e-03)(0.0125,6.118727e-04)(0.00625,1.596273e-04)};
\addplot[color=red,mark=*] coordinates {(0.1,2.822668e-02)(0.05,6.638873e-03)(0.025,1.687093e-03)(0.0125,4.080693e-04)(0.00625,1.021256e-04)};
\addplot[color=green,mark=*] coordinates {(0.1,3.805593e-02)(0.05,9.769536e-03)(0.025,2.663895e-03)(0.0125,6.185824e-04)(0.00625,1.603677e-04)};
\logLogSlopeTriangle{0.55}{0.2}{0.08}{2}{blue};
\legend{$L^2$,$H^1$,$L^\infty$}
\end{loglogaxis}
\end{tikzpicture}\\[-2pt]
{\footnotesize (b) Poisson, pure Neumann}
\end{minipage}\hfill
\caption{(a) Bean-domain $\Omega_h$ (dots; interface in red; the two boundary-isolated points of this coarse grid, $h=1/7$, in blue). (b) Relative $L^2$ (blue), $H^1$ (red), $L^\infty$ (green) errors vs $h$; fitted slopes $1.99$, $2.02$, $1.98$.}
\label{fig:results}
\end{figure}

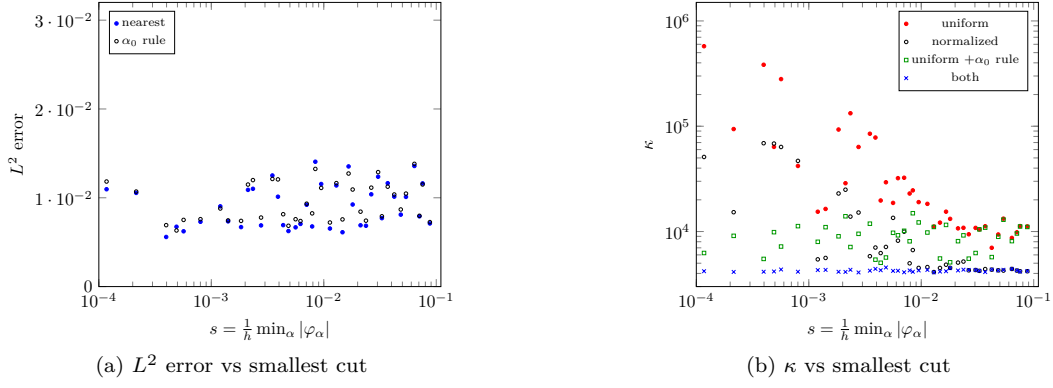
\begin{figure}[t]
\centering
\begin{minipage}[c]{0.49\linewidth}\centering
\begin{tikzpicture}[thick,scale=0.66, every node/.style={scale=1.0}]
\begin{semilogxaxis}[xlabel={$s=\frac{1}{h}\min_{\alpha}|\varphi_\alpha|$},ylabel=$L^2$ error,ymin=0,ymax=3.2e-2,xmin=1e-4,xmax=1.1e-1,legend style={at={(0.03,0.97)},anchor=north west,font=\scriptsize},scaled y ticks=false,ylabel near ticks]
\addplot[only marks,mark=*,mark size=1pt,color=blue] coordinates {(1.171e-04,1.096e-02)(2.148e-04,1.055e-02)(3.969e-04,5.587e-03)(4.903e-04,6.740e-03)(5.669e-04,6.225e-03)(8.008e-04,7.296e-03)(1.202e-03,9.035e-03)(1.410e-03,7.345e-03)(1.839e-03,6.699e-03)(2.119e-03,1.089e-02)(2.345e-03,1.101e-02)(2.769e-03,6.889e-03)(3.493e-03,1.251e-02)(3.914e-03,1.012e-02)(4.356e-03,6.920e-03)(4.851e-03,6.253e-03)(5.615e-03,6.671e-03)(6.206e-03,7.057e-03)(7.044e-03,9.214e-03)(7.885e-03,6.772e-03)(8.412e-03,1.406e-02)(9.484e-03,1.155e-02)(1.134e-02,6.539e-03)(1.294e-02,1.139e-02)(1.467e-02,6.118e-03)(1.661e-02,1.353e-02)(1.806e-02,9.242e-03)(2.122e-02,6.901e-03)(2.372e-02,6.848e-03)(2.644e-02,1.039e-02)(3.042e-02,1.237e-02)(3.287e-02,7.705e-03)(3.703e-02,1.164e-02)(4.233e-02,1.013e-02)(4.837e-02,8.102e-03)(5.380e-02,1.010e-02)(6.392e-02,1.358e-02)(7.086e-02,7.910e-03)(7.563e-02,1.162e-02)(8.776e-02,7.109e-03)};
\addplot[only marks,mark=o,mark size=1pt,color=black] coordinates {(1.171e-04,1.183e-02)(2.148e-04,1.070e-02)(3.969e-04,6.918e-03)(4.903e-04,6.320e-03)(5.669e-04,7.508e-03)(8.008e-04,7.592e-03)(1.202e-03,8.792e-03)(1.410e-03,7.451e-03)(1.839e-03,7.395e-03)(2.119e-03,1.149e-02)(2.345e-03,1.199e-02)(2.769e-03,7.773e-03)(3.493e-03,1.212e-02)(3.914e-03,1.207e-02)(4.356e-03,8.147e-03)(4.851e-03,6.847e-03)(5.615e-03,7.584e-03)(6.206e-03,7.381e-03)(7.044e-03,9.333e-03)(7.885e-03,8.245e-03)(8.412e-03,1.325e-02)(9.484e-03,1.111e-02)(1.134e-02,7.216e-03)(1.294e-02,1.165e-02)(1.467e-02,7.565e-03)(1.661e-02,1.273e-02)(1.806e-02,1.093e-02)(2.122e-02,8.435e-03)(2.372e-02,7.423e-03)(2.644e-02,1.114e-02)(3.042e-02,1.289e-02)(3.287e-02,7.915e-03)(3.703e-02,1.125e-02)(4.233e-02,1.037e-02)(4.837e-02,8.682e-03)(5.380e-02,1.047e-02)(6.392e-02,1.382e-02)(7.086e-02,7.971e-03)(7.563e-02,1.150e-02)(8.776e-02,7.262e-03)};
\legend{nearest,$\alpha_0$ rule}
\end{semilogxaxis}
\end{tikzpicture}\\[-2pt]
{\footnotesize (a) $L^2$ error vs smallest cut}
\end{minipage}\hfill
\begin{minipage}[c]{0.49\linewidth}\centering
\begin{tikzpicture}[thick,scale=0.66, every node/.style={scale=1.0}]
\begin{loglogaxis}[xlabel={$s=\frac{1}{h}\min_{\alpha}|\varphi_\alpha|$},ylabel=$\kappa$,xmin=1e-4,xmax=1.1e-1,ymin=3e3,ymax=1.5e6,legend style={at={(0.97,0.97)},anchor=north east,font=\scriptsize},ylabel near ticks]
\addplot[only marks,mark=*,mark size=1pt,color=red] coordinates {(1.171e-04,5.754e+05)(2.148e-04,9.402e+04)(3.969e-04,3.836e+05)(4.903e-04,6.365e+04)(5.669e-04,2.804e+05)(8.008e-04,4.195e+04)(1.202e-03,1.542e+04)(1.410e-03,1.639e+04)(1.839e-03,9.303e+04)(2.119e-03,2.875e+04)(2.345e-03,1.328e+05)(2.769e-03,6.355e+04)(3.493e-03,8.487e+04)(3.914e-03,7.811e+04)(4.356e-03,1.971e+04)(4.851e-03,2.937e+04)(5.615e-03,1.870e+04)(6.206e-03,3.213e+04)(7.044e-03,3.241e+04)(7.885e-03,2.296e+04)(8.412e-03,2.459e+04)(9.484e-03,1.903e+04)(1.134e-02,1.832e+04)(1.294e-02,1.115e+04)(1.467e-02,1.215e+04)(1.661e-02,1.543e+04)(1.806e-02,1.319e+04)(2.122e-02,1.072e+04)(2.372e-02,1.085e+04)(2.644e-02,9.418e+03)(3.042e-02,1.083e+04)(3.287e-02,1.048e+04)(3.703e-02,1.123e+04)(4.233e-02,7.006e+03)(4.837e-02,9.373e+03)(5.380e-02,1.325e+04)(6.392e-02,8.693e+03)(7.086e-02,9.835e+03)(7.563e-02,1.126e+04)(8.776e-02,1.118e+04)};
\addplot[only marks,mark=o,mark size=1pt,color=black] coordinates {(1.171e-04,5.106e+04)(2.148e-04,1.526e+04)(3.969e-04,6.898e+04)(4.903e-04,6.846e+04)(5.669e-04,6.355e+04)(8.008e-04,4.686e+04)(1.202e-03,5.442e+03)(1.410e-03,5.611e+03)(1.839e-03,2.298e+04)(2.119e-03,2.496e+04)(2.345e-03,1.386e+04)(2.769e-03,1.515e+04)(3.493e-03,5.830e+03)(3.914e-03,7.045e+03)(4.356e-03,6.241e+03)(4.851e-03,7.138e+03)(5.615e-03,1.348e+04)(6.206e-03,8.186e+03)(7.044e-03,9.967e+03)(7.885e-03,4.989e+03)(8.412e-03,6.668e+03)(9.484e-03,4.512e+03)(1.134e-02,4.596e+03)(1.294e-02,4.120e+03)(1.467e-02,4.507e+03)(1.661e-02,4.839e+03)(1.806e-02,4.518e+03)(2.122e-02,5.069e+03)(2.372e-02,5.185e+03)(2.644e-02,4.281e+03)(3.042e-02,4.312e+03)(3.287e-02,4.231e+03)(3.703e-02,4.400e+03)(4.233e-02,4.348e+03)(4.837e-02,4.280e+03)(5.380e-02,4.257e+03)(6.392e-02,4.412e+03)(7.086e-02,4.260e+03)(7.563e-02,4.189e+03)(8.776e-02,4.195e+03)};
\addplot[only marks,mark=square,mark size=1pt,color=green!60!black] coordinates {(1.171e-04,6.254e+03)(2.148e-04,9.100e+03)(3.969e-04,5.494e+03)(4.903e-04,9.861e+03)(5.669e-04,7.204e+03)(8.008e-04,1.125e+04)(1.202e-03,7.988e+03)(1.410e-03,1.100e+04)(1.839e-03,9.006e+03)(2.119e-03,1.396e+04)(2.345e-03,7.125e+03)(2.769e-03,9.512e+03)(3.493e-03,1.187e+04)(3.914e-03,5.406e+03)(4.356e-03,5.050e+03)(4.851e-03,5.665e+03)(5.615e-03,9.727e+03)(6.206e-03,9.206e+03)(7.044e-03,1.016e+04)(7.885e-03,8.076e+03)(8.412e-03,1.492e+04)(9.484e-03,1.220e+04)(1.134e-02,9.786e+03)(1.294e-02,1.105e+04)(1.467e-02,5.542e+03)(1.661e-02,1.156e+04)(1.806e-02,5.102e+03)(2.122e-02,8.097e+03)(2.372e-02,9.196e+03)(2.644e-02,5.524e+03)(3.042e-02,6.250e+03)(3.287e-02,1.048e+04)(3.703e-02,1.089e+04)(4.233e-02,5.715e+03)(4.837e-02,8.918e+03)(5.380e-02,1.296e+04)(6.392e-02,8.109e+03)(7.086e-02,9.578e+03)(7.563e-02,1.117e+04)(8.776e-02,1.105e+04)};
\addplot[only marks,mark=x,mark size=1.5pt,color=blue] coordinates {(1.171e-04,4.199e+03)(2.148e-04,4.133e+03)(3.969e-04,4.152e+03)(4.903e-04,4.174e+03)(5.669e-04,4.359e+03)(8.008e-04,4.153e+03)(1.202e-03,4.304e+03)(1.410e-03,4.321e+03)(1.839e-03,4.148e+03)(2.119e-03,4.135e+03)(2.345e-03,4.332e+03)(2.769e-03,4.088e+03)(3.493e-03,4.249e+03)(3.914e-03,4.419e+03)(4.356e-03,4.289e+03)(4.851e-03,4.556e+03)(5.615e-03,4.227e+03)(6.206e-03,4.245e+03)(7.044e-03,4.101e+03)(7.885e-03,4.271e+03)(8.412e-03,4.147e+03)(9.484e-03,4.269e+03)(1.134e-02,4.321e+03)(1.294e-02,4.094e+03)(1.467e-02,4.164e+03)(1.661e-02,4.195e+03)(1.806e-02,4.493e+03)(2.122e-02,4.267e+03)(2.372e-02,4.344e+03)(2.644e-02,4.311e+03)(3.042e-02,4.310e+03)(3.287e-02,4.237e+03)(3.703e-02,4.123e+03)(4.233e-02,4.355e+03)(4.837e-02,4.267e+03)(5.380e-02,4.330e+03)(6.392e-02,4.418e+03)(7.086e-02,4.269e+03)(7.563e-02,4.161e+03)(8.776e-02,4.212e+03)};
\legend{uniform,normalized,uniform $+\alpha_0$ rule,both}
\end{loglogaxis}
\end{tikzpicture}\\[-2pt]
{\footnotesize (b) $\kappa$ vs smallest cut}
\end{minipage}
\caption{Scalar Neumann, fixed $h=0.05$, many geometries (bean radius swept) so that the smallest coupled cell $s$ spans three decades. (a) The relative $L^2$ error varies by a factor $2.5$ (nearest interior node) or $2.2$ ($\alpha_0$ rule) with no trend. (b) The $2$-norm condition number for four ways of writing the boundary equations: scaled uniformly by a power of $h$ (full marks), normalized row by row (open marks), scaled uniformly with $x_{\alpha_0}$ kept away from the interface (squares), and both precautions together (crosses). Only the last is independent of the cut.}
\label{fig:cond}
\end{figure}

The code is freely available\footnote{\url{https://github.com/PhiFEM/publication_PhiFD_Neumann}}.
The test case uses the bean-shaped domain $\Omega=\{\varphi<0\}$ given in Figure~\ref{fig:results}(a), a slightly non-convex bent ellipse $\varphi(x,y)=0.8\,x^2+(y+1.4\,x^2)^2-R^2$, $R=0.5$ (the term $y+1.4\,x^2$ bends it along $y=-1.4\,x^2$).
We solve the scalar problem $-\Delta u+u=f$ with the Neumann condition $\partial_n u=g$ on the entire boundary $\partial\Omega$ 
and the manufactured solution $u=\sin(Kx)\cos(Ky)$, $K=\pi/0.8$; the right-hand side of the boundary rows of \eqref{eq:sys} is assembled from the data $G$ alone, as a user would. Figure~\ref{fig:results}(b) reports the relative errors in the discrete $L^2$, $H^1$, and $L^\infty$ norms: we obtain the optimal second-order convergence in all three norms. The scheme is also robust to weakly cut cells (Figure~\ref{fig:cond}(a)): when sweeping the geometry at a fixed mesh size such that the size $s=\frac{1}{h}\min_{\alpha\in \Omega_h}|\varphi_\alpha|$ of the smallest coupled cell spans three decades, the error varies by a factor $2.2$ with no trend (fitted slope $+0.03$ in $s$; factor $2.5$ if $x_{\alpha_0}$ is simply the nearest interior node); over the same sweep the cut-cell scheme of \cite{arias2018poisson} varies by $3.2$ in error and by $13.8$ in condition number, against $1.11$ below for $\varphi$-FD. The condition number behaves as $O(h^{-2})$ and, with the two precautions of Section~\ref{sec:scheme}, is insensitive to the smallest cut cell as well. Figure~\ref{fig:cond}(b) separates the two effects. Scaling all the boundary equations by the same power of $h$, chosen from the generic magnitude $|\varphi|\sim h$ so that they balance the Laplacian rows, makes $\kappa$ vary by a factor $82$ over the sweep, with a fitted slope $-0.54$ in $s$ (full marks). Normalizing each of them by its own sup-norm (times $h^{-2}$) instead --- a row scaling, so the errors of Figure~\ref{fig:cond}(a) are unchanged up to round-off --- reduces this to a factor $17$ and a slope $-0.43$ (open marks). Choosing $x_{\alpha_0}$ away from the interface then removes what remains: $\kappa$ lies between $4.09\times10^{3}$ and $4.56\times10^{3}$, a factor $1.11$ with a fitted slope $+0.00$ (crosses), for an unchanged error (Figure~\ref{fig:cond}(a)). The two precautions are complementary rather than redundant: the $\alpha_0$ rule alone, with the uniform scaling, already removes the trend (slope $+0.02$) but leaves a factor $2.95$ of scatter (squares), which the normalization then reduces to $1.11$; the rule removes the degeneracy, the normalization the uneven weighting. The growth reported by the first variant is therefore an artifact of how the boundary equations are written, not a property of the scheme.

\section{Conclusion}
We have extended $\varphi$-FD, an immersed finite-difference scheme on Cartesian grids, to Neumann boundary conditions, which are imposed directly through the level-set function. On a non-convex test domain, it achieves second-order convergence in the discrete $L^2$, $H^1$, and $L^\infty$ norms. The $O(h^2)$ consistency established above and the numerically observed $O(h^{-2})$ conditioning are consistent with this second-order convergence, but a full stability proof is left to future work; a variational reformulation with a consistent weak form of the boundary condition, in the spirit of $\varphi$-FEM, would make coercivity immediate, at the price of changing the scheme. Two elementary precautions --- normalizing each boundary equation individually and keeping the coupled interior node away from the interface --- make this conditioning independent of the smallest cut cell as well, so that no ghost-penalty stabilization is needed here, unlike in the Dirichlet case~\cite{phiFD}.

\FloatBarrier
\bibliographystyle{abbrv}
\bibliography{biblio}

\end{document}